\documentclass[12pt,leqno]{article}
 
\def\w{\dot{w}}

\def\sigmad{\dot{\sigma}} 
\def\int{\mathbb{Z}}

\def\OO{{\cal O}}

\def\SO{{\rm SO}}
\def\HSpin{{\rm HSpin}}
\def\Spin{{\rm Spin}}
\def\Sp{{\rm Sp}}
\def\PSO{{\rm PSO}}
\def\PSp{{\rm PSp}}

\def\proof{{\bf Proof. }}

\def\pf{\proof}
\usepackage{times}
\usepackage{graphics}
\usepackage{amssymb}
\usepackage{amscd}
\usepackage{amsmath}
\input xy 
\xyoption{all}
\title{Quotients for sheets of  conjugacy classes}
\newtheorem{theorem}{Theorem}[section]
\newtheorem{lemma}[theorem]{Lemma}

\newtheorem{remark}[theorem]{Remark}

\author{Giovanna Carnovale, Francesco Esposito\\
Dipartimento di Matematica ``Tullio Levi-Civita''\\
Torre Archimede - via Trieste 63 - 35121 Padova - Italy\\
email: carnoval@math.unipd.it, esposito@math.unipd.it }

\date{}

\begin{document}
\maketitle
\begin{abstract}
We provide a description of the orbit space of a sheet $S$  for the conjugation action of a complex simple simply connected algebraic group $G$. This is obtained by means of a bijection between $S/G$ and the quotient of  a shifted torus modulo the action of a subgroup of the Weyl group and it is the group analogue of a result due to Borho and Kraft. We also describe the normalisation of the categorical quotient $\overline{S}//G$ for arbitrary simple $G$ and give a necessary and sufficient condition for $\overline{S}//G$ to be normal in analogy to results of Borho, Kraft and Richardson. The example of $G_2$ is worked out in detail. 
\end{abstract}

\section{Introduction}
Sheets for the action of a connected algebraic group $G$ on a variety $X$ have their origin in the work of Kostant \cite{kostant}, who studied the union of regular orbits for the adjoint action on a semisimple Lie algebra, and in the work of Dixmier \cite{dixmier}. Sheets are the irreducible components of the level sets of $X$ consisting of points whose orbits have the same dimension. In a sense they provide a natural way to collect orbits in families in order to study properties of one orbit by looking at others in its family. For the adjoint action of a complex semisimple algebraic group $G$ on its Lie algebra they were deeply and systematically studied in \cite{bo,BK}. They were described as sets, their closure was well-understood, they were classified in terms of pairs consisting of a Levi subalgebra and suitable nilpotent orbit therein, and they were used to answer affirmatively to a question posed by Dixmier on the multiplicities in the module decomposition of the ring of regular functions of an adjoint orbit in $\mathfrak{sl}(n,{\mathbb C})$. If $G$ is classical then all sheets are smooth \cite{imhof,peterson}. The study of sheets in positive characteristic has appeared more recently in \cite{PS}.

In analogy to this construction, sheets of primitive ideals were introduced and studied by W. Borho and A. Joseph in \cite{BJ}, in order to describe the set of primitive ideals in a universal enveloping algebra as a countable union of maximal varieties. More recently, Losev in \cite{losev} has introduced the notion of birational sheet in a semisimple Lie algebra, he has shown that birational sheets form a partition of the Lie algebra and has applied this result in order to establish a version of the orbit method for semisimple Lie algebras. Sheets were also used in \cite{PT} in order to parametrise the set of $1$-dimensional representations of finite $W$-algebras, with some applications also to the theory of primitive ideals. Closures of sheets appear as associated varieties of affine vertex algebras, \cite{moreau-arakawa}.

In characterisitc zero, several results on quotients $S/G$ and $\overline{S}//G$, for a sheet $S$ were addressed: Katsylo has shown in \cite{katsylo} that $S/G$ has the structure of a quotient and is isomorphic to the quotient of an affine variety by the action of a finite group \cite{katsylo}; Borho has explicitly described the normalisation of $\overline{S}//G$ and Richardson, Broer, Douglass-R\"ohrle in \cite{richardson, broer,DR} have provided the list of the quotients $\overline{S}//G$ that are normal. 

Sheets for the conjugation action of $G$ on itself were studied in \cite{gio-espo} in the spirit of \cite{BK}. If $G$ is semisimple, they are parametrized in terms of pairs consisting of a Levi subgroup of parabolic subgroups and a suitable isolated conjugacy class therein. Here isolated means that the connected centraliser of the semisimple part of a representative is semisimple. An alternative parametrisation can be given in terms of triples consisting of a pseudo-Levi subgroup $M$ of $G$, a coset in $Z(M)/Z(M)^\circ$ and a suitable unipotent class in $M$. Pseudo-Levi subgroups are, in good characteristic, centralisers of semisimple elements and up to conjugation they are subroot subgroups whose root system has a base in the extended Dynkin diagram of $G$ \cite{mcninch-sommers}. It is also shown in \cite{gio} that sheets in $G$ are the irreducible components of the parts in Lusztig's partition introduced in \cite{lusztig}, whose construction is given in terms of Springer's correspondence. 

Also in the group case one wants to reach a good understanding of quotients of sheets. An analogue of Katsylo's theorem was obtained for sheets containing spherical conjugacy classes and all such sheets are shown to be smooth \cite{gio-espo2}. The proof in this case relies on specific properties of the intersection of spherical conjugacy classes with Bruhat double cosets, which do not hold for general classes. Therefore, a straightforward generalization to arbitrary sheets is not immediate.  Even in absence of a Katsylo type theorem, it is of interest to understand the orbit space  $S/G$. In this paper we address the question for $G$ simple provided $G$ is simply connected if the root system is of type $C$ or $D$. We give a bijection between the orbit space $S/G$ and a quotient of a shifted torus of the form $Z(M)^\circ s$ by the action of a subgroup $W(S)$ of the Weyl group,  giving a group analogue of  \cite[Theorem 3.6]{kraft},\cite[Satz 5.6]{bo}.  In most cases the group $W(S)$ does not depend on the unipotent part of the triple corresponding to the given sheet although it may depend on the isogeny type of $G$. This is one of the difficulties when passing from the Lie algebra case to the group case. 
The restriction on $G$ needed for the bijection depends on the symmetry of the extended Dynkin diagram in this case: type $C$ and $D$ are the only two situations in which two distinct subsets of the extended Dynkin diagram can be equivalent even if they are not of type $A$. We illustrate by an example in $\HSpin_{10}({\mathbb C})$ that the restriction we put is necessary in order to have injectivity so our theorem is somehow optimal. 

We also address some questions related to the categorial quotient $\overline{S}//G$, for a sheet in $G$. We obtain group analogues of the description of the normalisation of $\overline{S}//G$ from \cite{bo} and of a necessary and sufficient condition on $\overline{S}//G$ to be normal from \cite{richardson}. Finally we apply our results to compute the quotients $S/G$ of all sheets in $G$ of type $G_2$ and verify which of the quotients $\overline{S}//G$ are normal. This example will serve as a toy example for a forthcoming paper in which we will list all normal quotients for $G$ simple. 

\section{Basic notions}
In this paper $G$ is a complex {\em simple} algebraic group with maximal torus $T$, root system $\Phi$, weight lattice $\Lambda$, set of simple roots $\Delta=\{\alpha_1,\,\ldots,\,\alpha_\ell\}$, Weyl group $W=N(T)/T$ and corresponding Borel subgroup $B$. The numbering of simple roots is as in \cite{bourbaki}. Root subgroups are denoted by $X_\alpha$ for $\alpha\in\Phi$ and their elements have the form $x_\alpha(\xi)$ for $\xi\in{\mathbb C}$. Let $-\alpha_0$ be the highest root and let $\tilde{\Delta}=\Delta\cup\{\alpha_0\}$. The centraliser of an element $h$ in a closed group $H\leq G$ will be denoted by $H^h$ and the identity component of $H$ will be indicated by $H^\circ$. 
If $\Pi\subset \tilde{\Delta}$ we set
\begin{equation*}
G_\Pi:=\langle T,\,X_{\pm\alpha}~|~\alpha\in\Pi\rangle.
\end{equation*}
Conjugates of such groups are called pseudo-Levi subgroups.
We recall from \cite[\S 6]{mcninch-sommers}  that if $s\in T$ then its connected centraliser $G^{s\circ}$ is conjugated to $G_\Pi$ for some $\Pi$ by means of an element in $N(T)$.  By \cite[2.2]{hu-cc} we have $G^s=\langle G^{s\circ}, N(T)^s\rangle$. $W_\Pi$ indicates the subgroup of $W$ generated by the simple reflections with respect to roots in $\Pi$ and it is the Weyl group of $G_\Pi$. 

We realize the groups $\Sp_{2\ell}({\mathbb C})$,  $\SO_{2\ell}({\mathbb C})$ and $\SO_{2\ell+1}({\mathbb C})$, respectively, as the groups of matrices of determinant $1$ preserving the bilinear forms: $\left(\begin{smallmatrix} 0&I_\ell
\\ -I_\ell&0\end{smallmatrix}\right)$,  $\left(\begin{smallmatrix} 0&I_\ell
\\ I_\ell&0\end{smallmatrix}\right)$ and  $\left(\begin{smallmatrix}
1\\
&&I_\ell\\
&I_\ell\\
\end{smallmatrix}\right)$, respectively.

If $G$ acts on a variety $X$, the action of $g\in G$ on $x\in X$ will be indicated by $(g,x)\mapsto g\cdot x$.  If $X=G$ with adjoint action we thus have $g\cdot h=ghg^{-1}$. 
For $n\geq 0$ we shall denote by $X_{(n)}$ the union of orbits of dimension $n$. The nonempty sets $X_{(n)}$ are locally closed and a sheet $S$ for the action of $G$ on $X$ is an irreducible component of any of these. For any $Y\subset X$ we set $Y^{reg}$ to be the set of points of $Y$ whose orbit has maximal dimension. We recall the parametrisation and description of sheets for the action of $G$ on itself by conjugation and provide the necessary background material. 

 A Jordan class  in $G$ is an equivalence class with respect to the equivalence relation:
 $g,\,h\in G$ with Jordan decomposition  $g=su$, $h=rv$ are equivalent if up to conjugation $G^{s\circ}=G^{r\circ}$, $r\in Z(G^{s\circ})^\circ s$ and $G^{s\circ}\cdot u=G^{s\circ}\cdot v$. 
 As a set, the Jordan class of $g=su$ is thus $J(su)=G\cdot((Z(G^{s\circ})^\circ s)^{reg}u)$ and it is contained in some $G_{(n)}$. Jordan classes are parametrised by $G$-conjugacy classes of triples $(M, Z(M)^\circ s, M\cdot u)$ where $M$ is a pseudo-Levi subgroup, $Z(M)^\circ s$ is a coset in $Z(M)/Z(M)^\circ$ such that $(Z(M)^\circ s) ^{reg}\subset Z(M)^{reg}$ and $M\cdot u$ is a unipotent conjugacy class in $M$. They are finitely many, locally closed, $G$-stable, smooth, see \cite[3.1]{lusztig-inventiones} and \cite[\S 4]{gio-espo} for further details. 

Every sheet $S\subset G$ contains a unique dense Jordan class, so sheets are parametrised by conjugacy classes of a subset of the triples above mentioned. 
More precisely, a Jordan class $J=J(su)$ is dense in a sheet if and only if it is not contained in $(\overline{J'})^{reg}$ for any Jordan class $J'$ different from $J$. We recall from \cite[Proposition 4.8]{gio-espo} that 
\begin{equation}
\overline{J(su)}^{reg}=\bigcup_{z\in Z(G^{s\circ})^\circ }G\cdot(s {\rm Ind}_{G^{s\circ}}^{G^{zs\circ}}(G^{s\circ}\cdot u)),
\end{equation}
where ${\rm Ind}_{G^{s\circ}}^{G^{zs\circ}}(G^{s\circ}\cdot u)$ is Lusztig-Spaltenstein's induction from the Levi subgroup $G^{s\circ}$ of a parabolic subgroup of $G^{zs\circ}$ of the class of $u$ in $G^{s\circ}$, see \cite{LS}. So, Jordan classes that are dense in a sheet correspond to triples where $u$ is a rigid orbit in $G^{s\circ}$, i.e., such that its class in $G^{s\circ}$ is not induced from a conjugacy class in a proper Levi subgroup of a parabolic subgroup of $G^{s\circ}$. 

A sheet consists of a single conjugacy class if and only if $\overline{S}=\overline{J(su)}=\overline{G\cdot su}$ where $u$ is rigid in $G^{s\circ}$ and $G^{s\circ}$ is semisimple, i.e., if and only if $s$ is isolated and $u$ is rigid in $G^{s\circ}$. Any sheet $S$ in $G$ is the image through the isogeny map $\pi$ of a sheet $S'$  in   the simply-connected cover $G_{sc}$ of $G$, where $S'$ is defined up to multiplication by an element in ${\rm Ker}(\pi)$. Also,  $Z(G^{\pi(s)\circ})=\pi(Z(G_{sc}^{s\circ}))$ and  $Z(G^{\pi(s)\circ})^\circ=\pi(Z(G_{sc}^{s\circ})^\circ)=Z(G_{sc}^{s\circ}) ^\circ{\rm Ker}(\pi)$.

\section{A parametrization of orbits in a sheet}

In this section we parametrize the set $S/G$ of conjugacy classes in a given sheet.
Let $S=\overline{J(su)}^{reg}$ with $s\in T$ and $u\in U\cap G^{s\circ}$. Let $Z=Z(G^{s\circ})$ and $L=C_G(Z^\circ)$. The latter is always a Levi subgroup of a parabolic subgroup of $G$, \cite[Proposition 8.4.5, Theorem 13.4.2]{springer} and if $\Psi_s$ is the root system of $G^{s\circ}$  with respect to $T$, then $L$ has root system $\Psi:={\mathbb Q}\Psi_s\cap \Phi$.

Let \begin{equation}\label{eq:WS}W(S)=\{w\in W~|~w(Z^\circ s)=Z^\circ s\}.\end{equation} We recall that $C_G(Z(G^{s\circ})^\circ s)^\circ=G^{s\circ}$. Thus, for any lift $\w$ of  $w\in W(S)$ we have $\w\cdot G^{s\circ}=G^{s\circ}$, so $\w\cdot Z^\circ=Z^\circ$ and therefore $\w\cdot L=L$. Thus, any $w\in W(S)$ determines an automorphism of $\Psi_s$ and $\Psi$. Let $\OO=G^{s\circ}\cdot u$. We set:
 \begin{equation}\label{eq:WSu}W(S)^u=\{w\in W(S)~|~\w\cdot\OO=\OO\}.\end{equation} 
The definition is independent of the choice of the representative of each $w$ because $T\subset L$. 

\begin{lemma}\label{lem:WS1}Let $\Psi_s$ be the root system of $G^{s\circ}$ with respect to $T$, with basis $\Pi\subset \Delta\cup\{-\alpha_0\}$. Let $W_\Pi$ be the Weyl group of $G^{s\circ}$ and let $W^\Pi=\{w\in W~|~w\Pi=\Pi\}$. Then
$$W(S)=W_\Pi\rtimes (W^\Pi)_{Z^\circ s}=\{w\in W_\Pi W^\Pi~|~wZ^\circ s=Z^\circ s\}.$$
In particular, if $G^{s\circ}$ is a Levi subgroup of a parabolic subgroup of $G$, then $W(S)=W_\Pi\rtimes W^\Pi=N_W(W_\Pi)$ and it is independent of the isogeny class of $G$. 
\end{lemma}
\pf Let $W_X$ denote the stabilizer of $X$ in $W$ for $X=Z^\circ s, G^{s\circ}, Z, Z^\circ$. We have the following chain of inclusions:
$$W(S)=W_{Z^\circ s}\leq W_{G^{s\circ}}\leq W_Z\leq W_{Z^\circ}.$$
We claim that $W_{G^{s\circ}}=W_\Pi\rtimes W^\Pi$. Indeed, $W_\Pi W^\Pi\leq W_{G^{s\circ}}$ is immediate and if  $w\in W_{G^{s\circ}}$ then $w\Psi_s=\Psi_s$ and $w\Pi$ is a basis for $\Psi_s$. Hence, there is  some $\sigma\in W_\Pi$ such that $\sigma w\in W^\Pi$. By construction $W^\Pi$ normalises $W_\Pi$.  The elements of $W_{G^{s\circ}}$ permute the connected components of $Z=Z(G^{s\circ})$ and $W_{Z^\circ s}$ is precisely the stabilizer of $Z^\circ s$ in there. Since the elements of $W_\Pi$ fix the elements of $Z(G^{s\circ})$ pointwise, they stabilize $Z^\circ s$, whence the statement. The last statement follows from the equality $W_\Pi\ltimes W^\Pi=N_W(W_\Pi)$ in \cite[Corollary3]{howl} and  \cite[Lemma 33]{mcninch-sommers} because in this case $Z^\circ s= zZ^\circ$ for some $z\in Z(G)$, so  $W_{Z^\circ s}=W_{Z^\circ}$.\hfill$\Box$ 

\begin{remark}If $G^{s\circ}$ is not a Levi subgroup of a parabolic subgroup of $G$, then $W(S)$ might depend on the isogeny type of $G$. For instance, if  $\Phi$ is of type $C_5$ and $s={\rm diag}(-I_2,x,I_2,-I_2, x^{-1},I_2)\in \Sp_{10}({\mathbb C})$ for $x^2\neq1$, then: 
\begin{align*}
&\Pi=\{\alpha_0,\alpha_1,\alpha_4,\,\alpha_5\}\\
&Z=Z(G^{s\circ})=\{{\rm diag}(\epsilon I_2,y, \eta I_2,\epsilon I_2,y^{-1},\eta I_2),\;y\in{\mathbb C}^*, \epsilon^2=\eta^2=1\},\\
&Z^\circ s=\{{\rm diag}(-I_2,I_2,y,-I_2,I_2, y^{-1}),\;y\in{\mathbb C}^*\},
\end{align*}
and $W^\Pi=\langle s_{\alpha_1+\alpha_2+\alpha_3+\alpha_4} s_{\alpha_2+\alpha_3}\rangle$. Since $s_{\alpha_1+\alpha_2+\alpha_3+\alpha_4} s_{\alpha_2+\alpha_3}(Z^\circ s)=-Z^\circ s$ we have  $W(S)=W_\Pi$. However, if $\pi\colon \Sp_{10}({\mathbb C})\to \PSp_{10}({\mathbb C}) $ is the isogeny map, then
$W^\Pi$ preserves $\pi(Z^\circ s)$ so $W(\pi(S))=W_\Pi\rtimes W^\Pi$. Taking $u=1$ have an example in which also  $W(S)^u$ depends on the isogeny type.
\end{remark}

\begin{table}[ht]
\caption{Kernel of the isogeny map; $\Phi$ of type $B_\ell$,  $C_\ell$ or $D_\ell$}\label{tab:kernel}
\vskip0.5cm
\begin{tabular}{|c|c|c|p{6cm}|}
\hline  type  & parity of $\ell$ &  group   &  ${\rm Ker}\pi$ \\
\hline
$B_\ell$ & any & $\SO_{2\ell+1}({\mathbb C})$ & $\left\langle \alpha^{\vee}_\ell(-1)\right\rangle$\\
\hline
 $C_\ell$ & any  & $\PSp_{2\ell}({\mathbb C})$ & $\left\langle \displaystyle\prod_{j\text{ odd}}\alpha^{\vee}_j(-1)\right\rangle=\langle-I_{2\ell}\rangle$\\
 \hline 
 $D_{\ell}$ & even & $\PSO_{2\ell}({\mathbb C})$ &$ \left\langle \displaystyle\prod_{j\text{ odd}} \alpha^{\vee}_j(-1),\,\alpha^{\vee}_{\ell-1}(-1)\alpha^{\vee}_{\ell}(-1)\right\rangle$\\
\hline
 $D_{\ell}$ & odd &  $\PSO_{2\ell}({\mathbb C})$ &  $\left\langle \displaystyle \prod_{j\text{ odd}\leq \ell - 2}\alpha^{\vee}_j (-1)\alpha^{\vee}_{\ell -1}(i)\alpha^{\vee}_{\ell}(i^3)\right\rangle$\\
\hline
 $D_{\ell}$ & any &  $\SO_{2\ell}({\mathbb C})$ & $\left\langle  \alpha^{\vee}_{\ell-1}(-1)\alpha^{\vee}_{\ell}(-1)\right\rangle$\\
\hline
 $D_{\ell}$ & even &  $\HSpin_{2\ell}({\mathbb C})$ &  $\left\langle \displaystyle \prod_{j\text{ odd}}\alpha^{\vee}_j (-1))\right\rangle$\\
\hline
\end{tabular}
\end{table}

Next Lemma shows that in most cases $W(S)^u$ can be determined without any knowledge of $u$.

\begin{lemma}\label{WS}Suppose $G$ and $S=\overline{J(su)}^{reg}$ are {\bf not} in the following situation: 
\begin{center}
``$G$ is either $\PSp_{2\ell}({\mathbb C})$, $\HSpin_{2\ell}({\mathbb C})$, or $\PSO_{2\ell}({\mathbb C})$;\\
 $[G^{s\circ},\,G^{s\circ}]$ has two isomorphic simple factors $G_1$ and $G_2$ that are not of type $A$;\\ the components of $u$ in $G_1$ and $G_2$ do not correspond to the same partition.''
\end{center}
Then, $W(S)=W(S)^u$. \end{lemma}
\pf The element $u$ is rigid in $[G^{s\circ},G^{s\circ}]\leq G^{s\circ}$ and this happens if and only if each of its components in the corresponding simple factor of $[G^{s\circ},G^{s\circ}]$ is rigid.  Rigid unipotent elements in type $A$ are trivial \cite[Proposition 5.14]{spaltenstein}, therefore what matters are only the components of $u$ in the simple factors of type different from $A$. In addition, 
rigid unipotent classes are characteristic in simple groups, \cite[Lemma 3.9, Korollar 3.10]{bo}.  For all $\Phi$ different from $C$ and $D$, simple factors that are not of type A are never isomorphic. Therefore the statement certainly holds in all cases with a possible exception when: $\Phi$ is of type $C_\ell$ or $D_\ell$; $[G^{s\circ},G^{s\circ}]$ has two isomorphic factors  of type different from $A$; and the components of $u$ in those two factors, that are of type $C_m$ or $D_m$, respectively,  correspond to different partitions.

Let us assume that we are in this situation. Then, $W(S)=W(S)^u$ if and only if the elements of $W(S)$, acting as automorphisms of $\Psi_s$, do not interchange the two isomorphic factors in question. We have $2$ isogeny classes in type $C_{\ell}$, 3 in type $D_\ell$ for $\ell$ odd, and 4 (up to isomorphism) in type $D_\ell$ for $\ell$ even. 

If $\Phi$ is of type $C_\ell$ and $G=\Sp_{2\ell}({\mathbb C})$ up to a central factor  $s$ can be chosen to be of the form:
\begin{equation}\label{eq:esse}s={\rm diag}(I_{m}, t,-I_{m},I_m, t^{-1},-I_m)\end{equation} where $t$ is a diagonal matrix in ${\rm GL}_{\ell-2m}({\mathbb C})$ with eigenvalues different from $\pm1$. 
Then $\Pi$ is the union of $\{\alpha_0,\,\ldots,\,\alpha_{m-1}\}$, $\{\alpha_\ell,\alpha_{\ell-1},\ldots,\,\alpha_{\ell-m+1}\}$ and possibly other subsets of simple roots orthogonal to these.
Then $W^{\Pi}$ is the direct product of terms permuting isomorphic components of type $A$ with the subgroup generated by $\sigma=\prod_{j=1}^m s_{\alpha_j+\cdots+\alpha_{\ell-j}}$. 
In this case the elements of  $Z^\circ s$ are of the form ${\rm diag}(I_{m}, r,-I_{m},I_m, r^{-1},-I_m)$, where $r$ has the same shape as $t$ and $\sigma(Z^\circ s)=-Z^\circ s$.
Thus, $W^\Pi$ does not permute the two factors of type $C_m$ and $W(S)=W(S)^u$. 

If, instead, $G=\PSp_{2\ell}({\mathbb C})$ and the sheet is $\pi(S)$, we may take $J=J(\pi(su))$ where $s$ is as in \eqref{eq:esse}. Then, $\sigma$ preserves $\pi(Z^\circ s)$ and therefore $W(\pi(S))\not=W(\pi(S))^{\pi(u)}$.

Let now $\Phi$ be of type $D_\ell$ and  $G=Spin_{2\ell}({\mathbb C})$. With notation as in \cite{springer}, we may take
\begin{equation}\label{eq:esseD}s=\left(\prod_{j=1}^m\alpha_j^\vee(\epsilon^j)\right)\left(\prod_{i=m+1}^{l-m-1}\alpha_i^\vee(c_i)\right)\left(\prod_{b=2}^m \alpha_{\ell-b}^\vee(d^2\eta^{b})\right)\alpha_{\ell-1}^\vee(\eta d)\alpha_\ell(d)\end{equation}
%
with $\epsilon^2=\eta^2=1$, $\epsilon\neq\eta$, and $d,\,c_i\in{\mathbb C}^*$. \\
Here $\Pi$ is the union of $\{\alpha_0,\,\ldots,\,\alpha_{m-1}\}$, $\{\alpha_\ell,\alpha_{\ell-1},\ldots,\,\alpha_{\ell-m+1}\}$ and possibly other subsets of simple roots orthogonal to these. 
Then $W^{\Pi}$ is the direct product of terms permuting isomorphic components of type $A$ and $\langle \sigma\rangle$ where $\sigma=\prod_{j=1}^m s_{\alpha_j+\cdots+\alpha_{\ell-j+1}}$. 
The coset  $Z^\circ s=Z_{\epsilon,\eta}$ consists of elements of the same form as \eqref{eq:esseD} with constant value of $\epsilon$ and $\eta$, and $Z^\circ=Z_{1,1}$ consists of the elements of similar shape with $\eta=\epsilon=1$. Then $\sigma(Z_{\epsilon,\eta})=Z_{\eta,\epsilon}$, hence $\sigma\not\in W(S)$, so $W(S)$ preserves the components of $\Psi_s$ of type $D$ and $W(S)=W(S)^u$. 

If $\ell=2q$ and $G=\HSpin_{2\ell}({\mathbb C})$ and $\pi\colon \Spin_{2\ell}({\mathbb C})\to \HSpin_{2\ell}({\mathbb C})$ is the isogeny map we see from Table \ref{tab:kernel} that ${\rm Ker}(\pi)$ is generated by an element $k$ such that $kZ_{\epsilon,\eta}=Z_{-\epsilon,\eta}$, so  $\sigma$ as above preserves $\pi(Z^\circ s)$  whereas it does not preserve the conjugacy class of $\pi(u)$. Therefore $\sigma\in W(\pi(S))\neq W(\pi(S))^u$. 

If $G=\SO_{2\ell}({\mathbb C})$ and $\pi\colon \Spin_{2\ell}({\mathbb C})\to \SO_{2\ell}({\mathbb C})$ is the isogeny map, then ${\rm Ker}(\pi)$ is generated by an element $k$ such that $kZ_{\epsilon,\eta}=Z_{\epsilon,\eta}$. In this case $\sigma$ does not preserve $\pi(Z^\circ s)$, whence $\sigma\not\in W(\pi(S))= W(\pi(S))^u$. 

If $G=\PSO_{2\ell}({\mathbb C})$ and $\pi\colon \Spin_{2\ell}({\mathbb C})\to \PSO_{2\ell}({\mathbb C})$, then by the discussion of the previous isogeny types 
we see that $\sigma(Z_{\epsilon,\eta})\subset {\rm Ker}(\pi)Z_{\epsilon,\eta}$, so  $\sigma$  preserves $\pi(Z^\circ s)$ whence $\sigma\in W(\pi(S))\neq W(\pi(S))^u$. \hfill$\Box$
\bigskip

Following \cite[\S 5]{bo} and according to \cite[Proposition 4.7]{gio-espo} we define the map 
$$\begin{array}{rl}
\theta\colon Z^\circ s&\to S/G\\
zs&\mapsto {\rm Ind}_L^G(L\cdot szu)\end{array}$$
where $L=C_G(Z(G^{s\circ})^\circ)$. 

\begin{lemma}\label{invariant}With the above notation, $\theta(zs)=\theta(w\cdot(zs))$ for every $w\in W(S)^u$.
\end{lemma}
\pf Let us observe that, since $z\in Z(L)$ and $G^{s\circ}\subset L$ there holds $L^{zs\circ}=G^{s\circ}$. In particular, $G^{s\circ}$ is a Levi subgroup of a parabolic subgroup of $G^{zs\circ}$. Let $U_P$ be the unipotent radical of a parabolic subgroup of $G$ with Levi factor $L$ and let $\w$ be a representative of $w$ in $N_G(T)$. By \cite[Proposition 4.6]{gio-espo} we have
$$
\begin{array}{rl}
{\rm Ind}_L^G(L\cdot (w\cdot zs)u)&=G\cdot(w\cdot(zs) u U_P)^{reg}\\
&=G\cdot(zs (\w^{-1}\cdot u) U_{\w^{-1}\cdot P})^{reg}\\
&={\rm Ind}_L^G(L\cdot (zs (\w^{-1}\cdot u)))\\
&=G\cdot \left(zs\,{\rm Ind}_{G^{s\circ}}^{G^{zs\circ}}(\w^{-1}\cdot(G^{s\circ}\cdot u))\right)\\
&=G\cdot \left( zs\,{\rm Ind}_{G^{s\circ}}^{G^{zs\circ}}(G^{s\circ}\cdot u)\right)\\
&={\rm Ind}_L^G(L\cdot (zsu))
\end{array}
$$
where we have used that $L=\w\cdot L$ for every $w\in W(S)^u\leq W(S)$  and independence of the choice of the parabolic subgroup with Levi factor $L$, \cite[Proposition 4.5]{gio-espo}. \hfill$\Box$

\smallskip

\begin{remark}The requirement that $w$ lies in $W(S)^u$ rather than in $W(S)$ is necessary. For instance, we consider $G=\PSp_{2\ell}({\mathbb C})$ with $\ell=2m+1$ and $s$ the class of ${\rm diag}(I_m,\lambda,-I_m,I_m,\lambda^{-1},-I_m)$ with $\lambda^4\neq1$ and $u$ rigid with non-trivial component only in the  subgroup $H=\langle X_{\pm\alpha_j},\,j=0,\,\ldots m-1\rangle$ of $G^{s\circ}$. The element $\sigma=\prod_{j=1}^ms_{\alpha_j+\cdots+\alpha_{\ell-j}}$ lies in $W(S)\setminus W(S)^u$. Taking $\theta(s)$ we have 
$${\rm Ind}_L^G(L\cdot su)=G\cdot su$$ whereas 
$${\rm Ind}_L^G(L\cdot w(s)u)={\rm Ind}_L^G(L\cdot s (\w\cdot u))=G\cdot( s(\w\cdot u)),$$ where $\w$ is any representative of $w$ in $N_G(T)$.
These classes would coincide only if $u$ and $\w\cdot u$ were conjugate in $G^s$. They are not conjugate in $G^{s\circ}$ because they lie in different simple components. Moreover, $G^s$ is generated by $G^{s\circ}$ and the lifts of elements in the centraliser $W^s$ of $s$ in $W$ \cite[2.2]{hu-cc}, which is contained in $W(S)$. Since $\lambda^4\neq1$ we see that the elements of $W^s$ cannot interchange the two components of type $C_m$ in  $G^{s\circ}$. Hence, $$\theta(s)={\rm Ind}_L^G(L\cdot su)\neq{\rm Ind}_L^G(L\cdot w(s)u)=\theta(w(s)).$$
\end{remark}

In analogy with the Lie algebra case we formulate the following theorem. The proof follows the lines of \cite[Satz 5.6]{bo} but a more detailed analysis is necessary because the naive generalization of statement \cite[Lemma 5.4]{bo} from Levi subalgebras in a Levi subalgebra to Levi subgroups in a pseudo-Levi subgroup does not hold. 

\begin{theorem}\label{bijection}Assume $G$ is simple and different from  $\PSO_{2\ell}({\mathbb C})$, $\HSpin_{2\ell}({\mathbb C})$ and $\PSp_{2\ell}({\mathbb C})$, $\ell\geq5$. Let $S=\overline{J(su)}^{reg}$, with $s\in T$,  $Z=Z(G^{s\circ})$ and let $W(S)$  be as in \eqref{eq:WS}.  The map $\theta$ induces a bijection $\overline{\theta}$ between $Z^\circ s/W(S)$ and $S/G$. 
\end{theorem}
\pf Recall that under our assumptions Lemma \ref{WS} gives $W(S)=W(S)^u$. By Lemma \ref{invariant}, $\theta$ induces a well-defined map $\overline{\theta}\colon Z^\circ s/W(S)\to S/G$. It is surjective by \cite[Proposition 4.7]{gio-espo}. We prove injectivity.

Let us assume that $\theta(zs)=\theta(z's)$ for some $z,z'\in Z^\circ$. By construction, $Z^\circ\subset T$.
By \cite[Proposition 4.5]{gio-espo} we have
$$G\cdot \left(zs\left({\rm Ind}_{G^{s\circ}}^{G^{zs\circ}}(G^{s\circ}\cdot u)\right)\right)=G\cdot \left(z's\left({\rm Ind}_{G^{s\circ}}^{G^{z's\circ}}(G^{s\circ} \cdot u)\right)\right).$$
This implies that $z's=\sigma\cdot (zs)$ for some $\sigma\in W$. Let $\sigmad\in N(T)$ be a representative of $\sigma$. Then  
\begin{align*}\theta(zs)=\theta(z's)&=G\cdot \left((\sigma\cdot zs)({\rm Ind}_{G^{s\circ}}^{G^{z's\circ}}(G^{s\circ} \cdot u))\right)\\
&=G\cdot \left(zs\left({\rm Ind}_{\sigmad^{-1}\cdot (G^{s\circ})}^{\sigmad^{-1}\cdot G^{z's\circ}}\left(\sigmad^{-1}\cdot (G^{s\circ}\cdot u)\right)\right)\right)\\
&=G\cdot \left(zs\left({\rm Ind}_{\sigmad^{-1}\cdot (G^{s\circ})}^{G^{zs\circ}}\left(\sigmad^{-1}\cdot (G^{s\circ} \cdot u)\right)\right)\right).\end{align*}
Since the unipotent parts of $\theta(zs)$ and $\theta(z's)$ coincide, for some $x\in G^{zs}$ we  have 
$$x\cdot ({\rm Ind}_{G^{s\circ}}^{G^{zs\circ}}(G^{s\circ} \cdot u))={\rm Ind}_{\sigmad^{-1}\cdot (G^{s\circ})}^{G^{zs\circ}}\left(\sigmad^{-1}\cdot (G^{s\circ}\cdot u)\right).$$
The element $x$ may be written as $\w g$ for some $\w\in N(T)\cap G^{zs}$ and some $g\in G^{zs\circ}$ \cite[\S 2.2]{hu-cc}. Hence, 
\begin{align*}
{\rm Ind}_{G^{s\circ}}^{G^{zs\circ}}(G^{s\circ} \cdot u)&=\w^{-1}\cdot\left({\rm Ind}_{\sigmad^{-1}\cdot(G^{s\circ})}^{G^{zs\circ}}\left(\sigmad^{-1}\cdot(G^{s\circ} \cdot u)\right)\right)\\
&={\rm Ind}_{\w^{-1}\sigmad^{-1}\cdot (G^{s\circ})}^{G^{zs\circ}}\left((\w^{-1}\sigmad^{-1})\cdot(G^{s\circ} \cdot u)\right).
\end{align*}
Let us put
$$M:=G^{zs\circ}=\langle T, X_{\alpha},\alpha\in\Phi_M\rangle,\quad L_1:=G^{s\circ}= \langle T, X_{\alpha},\alpha\in\Psi\rangle$$ 
with $\Phi_M=\bigcup_{j=1}^l\Phi_j$ and $\Psi=\bigcup_{i=1}^m\Psi_i$ the decompositions in irreducible root subsystems. We recall that $L_1$ and $L_2:=(\w^{-1}\sigmad^{-1})\cdot L_1$ are Levi subgroups of some parabolic subgroups of $M$. We claim that if $L_1$ and $L_2$ are conjugate in $M$, then $zs$ and $z's$ are $W(S)$-conjugate. Indeed, under this assumption, since $L_1$ and $L_2$ contain $T$, there is $\dot{\tau}\in N_M(T)$ such that $L_1=\dot{\tau}\cdot L_2=\dot{\tau} \w^{-1} \sigmad^{-1}\cdot L_1$, so $\tau w^{-1}\sigma^{-1}(Z^\circ)=Z^\circ$. Then, $\tau w^{-1}\sigma^{-1}(z's)=zs$ and therefore 
$$\tau w^{-1}\sigma^{-1} (Z^\circ s)=\tau w^{-1}\sigma^{-1} (Z^\circ z's)=Z^\circ zs=Z^\circ s.$$ Hence $zs$ and $z's$ are $W(S)$-conjugate. By Lemma \ref{WS}, we have the claim. 
We show that if $\Phi_M$ has at most one component different from type $A$, then $L_1$ is always conjugate to $L_2$ in $M$. We analyse two possibilities.

\bigskip
\noindent \underline{$\Phi_j$ is of type $A$ for every $j$.} In this case the same holds for $\Psi_i$ and $u=1$. We recall that in type $A$ induction from the trivial orbit in a Levi subgroup corresponding to a partition $\lambda$ yields the unipotent class corresponding to the dual partition \cite[7.1]{spaltenstein}. Hence, equivalence of the induced orbits in each simple factor $M_i$ of $M$ forces 
$\Phi_j\cap\Psi\cong\Phi_j\cap w^{-1}\sigma^{-1}\Psi$ for every $j$. Invoking \cite[Lemma 5.5]{bo}, in each component $M_i$ we deduce that $L_1$ and $L_2$ are $M$-conjugate. 

\noindent\underline{There is exactly one component in $\Phi_M$ which is not of type $A$.} We set it to be $\Phi_1$. Then, there is at most one $\Psi_j$, say $\Psi_1$, which is not of type $A$, and $\Psi_1\subset \Phi_1$. In this case, $w^{-1}\sigma^{-1}\Phi_1\subset \Psi_1$.  Equivalence of the induced orbits in each simple factor $M_j$ of $M$ forces 
$\Phi_j\cap\Psi\cong\Phi_j\cap w^{-1}\sigma^{-1}\Psi$ for every $j>1$. By exclusion, the same isomorphism holds for $j=1$. Invoking once more \cite[Lemma 5.5]{bo} for each simple component, we deduce that $L_1$ and $L_2$ are $M$-conjugate. 

Assume now that there are exactly two components of $\Phi_M$ which are not of type $A$. This situation can only occur if $\Phi$ is of type $B_\ell$ for $\ell\geq6$, $C_\ell$ for $\ell\ge 4$ or $D_\ell$ for $\ell\geq 8$ (we recall that $D_2=A_1\times A_1$ and $D_3=A_3$). By a case-by-case analysis we directly show that $\sigma$ can be taken in $W(S)$. 

If $G={\rm Sp}_{2\ell}({\mathbb C})$  we may assume that
\begin{equation*}
s={\rm diag}(I_m,t,-I_p,I_m,t^{-1},-I_p)
\end{equation*}
with $p,\,m\geq2$ and $t$ a diagonal matrix with eigenvalues different from $0$ and $\pm1$. Then $Z^\circ s$ consists of matrices in this form, so $zs$ and $z's$ are of the form
$zs={\rm diag}(I_m,h,-I_p,I_m,h^{-1},-I_p)$ and  $z's={\rm diag}(I_m,g,-I_p,I_m,g^{-1},-I_p)$, where $h$ and $g$ are invertible diagonal matrices. The elements $zs$ and $z's$ are conjugate in $G$ if and only if  ${\rm diag}(h,\,h^{-1})$ and ${\rm diag}(g,\,g^{-1})$ are conjugate in $G'={\rm Sp}_{2(\ell-p-m)}({\mathbb C})$. This is the case if and only if they are conjugate in the normaliser of the torus $T'=G'\cap T$. The natural embedding $G'\to G$ given by 
\begin{align*}
\left(\begin{matrix}A&B\\
C&D\end{matrix}\right)
&\mapsto
\left(\begin{matrix}
I_m\\
&A&&B\\
&&I_{p+m}\\
&C&&D\\
&&&&I_p\end{matrix}\right)
\end{align*}
gives an embedding of $N_{G'}(T')\leq N_G(T)$ whose image lies in $W(S)$. Hence, $zs$ and $z's$ are necessarily $W(S)$-conjugate. This concludes the proof of injectivity for $G={\rm Sp}_{2\ell}({\mathbb C})$.

If $G={\rm Spin}_{2\ell+1}({\mathbb C})$, then we may assume that
\begin{equation*}
s=\left(\prod_{j=1}^{m}\alpha_j^\vee((-1)^j)\right)\left(\prod_{b=m+1}^{\ell-p-1}\alpha_b^\vee(c_b)\right)\left(\prod_{q=1}^{p}\alpha_{\ell-q}^\vee(c^2)\right)\alpha_\ell^\vee(c)
\end{equation*}
where $m\geq4$, $p\geq 2$, $c,\,c_b\in{\mathbb C}^*$ are generic. 
Here $Z^\circ s$ consists of elements of the form
\begin{equation*}
\left(\prod_{j=1}^{m}\alpha_j^\vee((-1)^j)\right)\left(\prod_{b=m+1}^{\ell-p-1}\alpha_b^\vee(d_b)\right)\left(\prod_{q=1}^{p}\alpha_{\ell-q}^\vee(d^2)\right)\alpha_\ell^\vee(d)
\end{equation*}
with $d_b,\,d \in{\mathbb C}^*$. 
The reflection $s_{\alpha_1+\cdots+\alpha_\ell}=s_{\varepsilon_1}$ maps any $y\in Z^\circ s$ to $y\alpha_\ell^\vee(-1)\in Z(G)Z^\circ s=Z^\circ s$. 

Let us consider the natural isogeny $\pi\colon G\to G_{ad}=\SO_{2\ell+1}({\mathbb C})$.  Then
\begin{equation*}
\pi(s)={\rm diag}(1,-I_{m},t,I_p,-I_m,t^{-1},I_p)
\end{equation*}
where $t$ is a diagonal matrix with eigenvalues different from $0$ and $\pm1$. A similar calculation as in the case of $\Sp_{2\ell}({\mathbb C})$ shows that $\pi(zs)$  is conjugate to $\pi(z's)$ by an element $\sigma_1\in W(\pi(S))=W(\pi(S))^u$. 
Then, $\sigma_1(zs)=kz's$, where $k\in Z(G)$. If $k=1$, then we set $\sigma=\sigma_1$ whereas if $k=\alpha_\ell^\vee(-1)$ we set $\sigma=s_{\alpha_1+\cdots+\alpha_\ell}\sigma_1$. Then $\sigma(zs)=z's$ and $\sigma(Z^\circ s)=Z(G)Z^\circ s=Z^\circ s$. This concludes the proof for ${\rm Spin}_{2\ell+1}({\mathbb C})$ and $\SO_{2\ell+1}({\mathbb C})$.

If $G={\rm Spin}_{2\ell}({\mathbb C})$, up to multiplication by a central element we may assume that
\begin{equation*}
s=\left(\prod_{j=m+1}^{\ell-p-1}\alpha_j^\vee(c_j)\right)\left(\prod_{q=2}^{p}\alpha_{\ell-q}^\vee((-1)^qc^2)\right)\alpha_{\ell-1}^\vee(-c)\alpha_{\ell}^\vee(c)
\end{equation*}
where $m,\,p\geq4$, $c,\,c_j\in{\mathbb C}^*$ are generic. The elements in $Z^\circ s$ are of the form
\begin{equation*}
\left(\prod_{j=m+1}^{\ell-p-1}\alpha_j^\vee(d_j)\right)\left(\prod_{q=2}^{p}\alpha_{\ell-q}^\vee((-1)^qd^2)\right)\alpha_{\ell-1}^\vee(- d)\alpha_{\ell}^\vee(d)
\end{equation*}
with $d_j,\,d \in{\mathbb C}^*$. We argue as we did for type $B_\ell$, considering the isogeny 
$\pi\colon G\to \SO_{2\ell}({\mathbb C})$.The Weyl group element $s_{\alpha_\ell}s_{\alpha_{\ell-1}}$ maps any $y\in Z^\circ s$ to $y\alpha_{\ell-1}^\vee(-1)\alpha_\ell^\vee(-1)\in {\rm Ker(\pi)}Z^\circ s=Z^\circ s$. 
The group $\pi(Z^\circ s)$ consists of elements of the form 
\begin{equation*}
{\rm diag}(I_m, t,-I_p,I_m,t^{-1},-I_p) 
\end{equation*}
where $t$ is a diagonal matrix in ${\rm GL}_{2(\ell-m-p)}({\mathbb C})$. Two elements 
\begin{align*}\pi(zs)={\rm diag}(I_m, h,-I_p,I_m,h^{-1},-I_p),\\
\pi(z's)={\rm diag}(I_m, g,-I_p,I_m,g^{-1},-I_p)\end{align*} therein are $W$-conjugate if and only if 
${\rm diag}(1, h,1,h^{-1})$ and $(1,g,1,g^{-1})$ are conjugate by an element $\sigma_1$ of the Weyl group $W'$ of $G'={\rm SO}_{2(\ell-m-p+1)}({\mathbb C})$. More precisely, even if $h$ and $g$ may have eigenvalues equal to $1$, we may choose $\sigma_1$ in the subgroup of $W'$ that either fixes the first and the $(\ell-m-p+2)$-th eigenvalues or interchanges them. Considering the natural embedding of $G'$ into ${\rm SO}_{2\ell}({\mathbb C})$ in a similar fashion as we did for $\SO_{2\ell}({\mathbb C})$, we show that $\sigma_1\in W(\pi(S))$. This proves injectivity for $\SO_{2\ell}({\mathbb C})$. Arguing as we did for ${\rm Spin}_{2\ell+1}({\mathbb C})$ using $s_{\alpha_\ell}s_{\alpha_{\ell-1}}$ concludes the proof of injectivity for ${\rm Spin}_{2\ell}({\mathbb C})$. 
\hfill$\Box$

\bigskip

The translation isomorphism $Z^\circ s\to Z^\circ$ determines a $W(S)$-equivariant map where $Z^\circ$ is endowed with the action $w\bullet z=(w\cdot zs)s^{-1}$, which is in general not  an action by automorphisms on $Z^\circ$.  Hence, $S/G$ is in bijection with the quotient $Z^\circ/W(S)$ of the torus $Z^\circ$ where the quotient is with respect to the $\bullet$ action. 

\bigskip

\begin{remark}{\rm Injectivity of $\overline{\theta}$ does not  necessarily hold for the adjoint groups $G=\PSp_{2\ell}({\mathbb C})$, $\PSO_{2\ell}({\mathbb C})$ and for $G=\HSpin_{2\ell}({\mathbb C})$. We give an example for $G=\HSpin_{20}({\mathbb C})$, in which $W(S)=W(S)^u$ and $G^{s\circ}$ is a Levi subgroup of a parabolic subgroup of $G$. Let $\pi\colon \Spin_{20}({\mathbb C})\to G$ be the central isogeny with kernel $K$ as in Table \ref{tab:kernel}.
Let $u=1$ and $$s=\alpha^\vee_{1}(a)\alpha_2^\vee(a^2)\alpha_3^\vee(a^3)\alpha_4^\vee(b)\alpha_5^\vee(c)\alpha_6^\vee(d^{-2}e^2)\alpha_7^\vee(e)\alpha_8^\vee(d^2)\alpha_9^\vee(d)\alpha_{10}^\vee(-d)K$$
with 
$a,\,b,\,c,\,d,\,e\in{\mathbb C}^*$ sufficiently generic. 
Then, $G^{s\circ}$ is generated by $T$ and the root subgroups of the subsystem with basis indexed by the following subset of the extended Dynkin diagram:
$$\begin{array}{lll}
\bullet-&\bullet-\circ-\circ-\circ-\circ-\bullet-&\circ-\bullet\\
&|&|\\
&\circ&\bullet\\
\end{array}
$$
Here $Z^\circ$ is given by elements of  shape:
$$\alpha^\vee_{1}(a_1)\alpha_2^\vee(a_1^2)\alpha_3^\vee(a_1^3)\alpha_4^\vee(b_1)\alpha_5^\vee(c_1)\alpha_6^\vee(d_1^{-2}e_1^2)\alpha_7^\vee(e_1)\alpha_8^\vee(d_1^2)\alpha_9^\vee(d_1)\alpha_{10}^\vee(-d_1)K$$ with $a_1,\,b_1,\,c_1,\,d_1,\,e_1\in{\mathbb C}^*$. Let 
$$zs=\alpha_5^\vee(c)\alpha_6^\vee(d^2)\alpha_7^\vee(-d^2)\alpha_8^\vee(d^2)\alpha_9^\vee(d)\alpha_{10}^\vee(-d)K\in Z^\circ sK$$ obtained by setting $a_1=b_1=1$, $c_1=c$, $d_1=d$ and $e_1=-d^2$,  and 
$$z's=\alpha_5^\vee(-c)\alpha_6^\vee(d^2)\alpha_7^\vee(-d^2)\alpha_8^\vee(d^2)\alpha_9^\vee(d)\alpha_{10}^\vee(-d)K\in Z^\circ sK,$$ obtained by setting $a_1=b_1=1$, $c_1=-c$, $d_1=d$ and $e_1=-d^2$.
%
The subgroup $M:=G^{zs\circ}=G^{z's\circ}$ is generated by $T$ and the root subgroups of the subsystem with basis indexed by the following subset of the extended Dynkin diagram:
$$\begin{array}{lll}
\bullet-&\bullet-\bullet-\circ-\circ-\circ-\bullet-&\bullet-\bullet\\
&|&|\\
&\bullet&\bullet\\
\end{array}
$$
For $\sigma=\prod_{j=1}^4s_{\alpha_j+\cdots+\alpha_{10-j}}$ we have $\sigma\cdot zs=z's$. We claim that $zs$ and $z's$ are not $W(S)$-conjugate. Equivalently, we show that  $\sigma W^{zsK}\cap W(S)=\emptyset$, where $W^{szK}$ is the stabiliser of $zs$ in $W$. Let $\sigma w$ be an element lying in such an intersection. We observe that if $\sigma w\in W(S)$, then $\sigma w(G^{s\circ})=G^{s\circ}$ hence $\sigma w$ cannot interchange the component of type $3A_1$ with the component of type $A_2$ therein. Thus, it cannot interchange the two components of type $D_4$ in $M$.  However, by looking at the projection $\pi'$ onto $G/Z(G)=\PSO_{10}({\mathbb C})$,  we see that $zsZ(G)$ is the class of the matrix 
\begin{align*}
{\rm diag}(I_4,c,c^{-1}d^2,-I_4,I_4,d^{-2}c,c^{-1},-I_4)
\end{align*}
which cannot be centralized by a Weyl group element interchanging these two factors if $c$ and $d$ are sufficiently generic. A fortiori, this cannot happen  for the class $zsK$. Hence, $zs$ and $z's$ are not $W(S)$-conjugate.

Let now $M_1$ and $M_2$ be the simple factors of $M$ corresponding respectively to the roots $\{\alpha_j,\,0\leq j\leq 3\}$, and  $\{\alpha_k,\,7\leq k\leq 10\}$, let $L_1=M_1\cap G^{s\circ}$ and $L_2=M_2\cap G^{s\circ}$.
Then, 
$$\theta(zs)={\rm Ind}_L^G(L\cdot zs)=G\cdot \left(zs({\rm Ind}_{G^{s\circ}}^M(1))\right)=G\cdot (zs({\rm Ind}_{L_1}^{M_1}(1))({\rm Ind}_{L_2}^{M_2}(1)))$$ and 
$$\theta(z's)={\rm Ind}_L^G(L\cdot z's)=G\cdot \left(z's({\rm Ind}_{G^{s\circ}}^M(1)\right)=G\cdot (z's({\rm Ind}_{L_1}^{M_1}(1))({\rm Ind}_{L_2}^{M_2}(1))).$$
Since $\sigma(zs)=z's$ we have, for some representative $\sigmad\in N(T)$:
$$\begin{array}{rl}
\theta(z's)&=G\cdot \left(zs ({\rm Ind}_{\sigmad^{-1} \cdot L_1}^{\sigmad^{-1}\cdot M_1}(1))({\rm Ind}_{\sigmad^{-1}\cdot L_2}^{\sigmad^{-1}\cdot M_2}(1)))\right)\\
&=G\cdot \left(zs ({\rm Ind}_{\sigmad^{-1} \cdot L_1}^{M_2}(1))({\rm Ind}_{\sigmad^{-1}\cdot L_2}^{M_1}(1)))\right).
\end{array}
$$
By \cite[Example 3.1]{moreau} we have 
${\rm Ind}_{\sigmad^{-1} \cdot L_1}^{M_2}(1)={\rm Ind}_{L_1}^{M_2}(1)$
and ${\rm Ind}_{\sigmad^{-1}\cdot L_2}^{M_1}(1)={\rm Ind}_{L_1}^{M_1}(1)$
so $\theta(zs)=\theta(z's)$.}\end{remark}

\begin{remark}The parametrisation in Theorem \ref{bijection} cannot be directly generalised to arbitrary Jordan classes. Indeed, if $u\in L$ is not rigid, then $L\cdot u$ is not necessarily characteristic and it may happen that for some external automorphism $\tau$ of $L$, the class $\tau(L\cdot u)$ differs from $L\cdot u$ even if they induce the same $G$-orbit. Then the map $\overline{\theta}$ is not necessarily injective.
\end{remark}

\section{The quotient $\overline{S}//G$}

In this section we discuss some properties of the categorical quotient $\overline{S}//G={\rm Spec}({\mathbb C}[\overline{S}])^G$ for $G$ simple in any isogeny class. 
Since  $\overline{S}//G$ parametrises only semisimple conjugacy classes it is enough to look at the so-called Dixmier sheets, i.e., the sheets containing a dense Jordan class consisting of semisimple elements. In addition, since every such Jordan class is dense in some sheet, studying the collection of $\overline{S}//G$ for $S$ a sheet in $G$ is the same as studying the collection of $\overline{J(s)}//G$ for $J(s)$ a semisimple Jordan class in $G$. 

The following Theorem is a group version of \cite[Satz 6.3]{bo}, \cite[Theorem 3.6(c)]{kraft} and \cite[Theorem A]{richardson}.
\begin{theorem}
Let $S=\overline{J(s)}^{reg}\subset G$. 
\begin{enumerate}
 \item  The normalisation of $\overline{S}//G$ is $Z(G^{s\circ})^\circ s/W(S)$. 
 \item The variety $\overline{S}//G$ is normal if and only if the natural map 
 \begin{equation}\label{eq:rho}\rho\colon {\mathbb C}[T]^W\to {\mathbb C}[Z(G^{s\circ})^\circ s]^{W(S)}\end{equation} induced from the restriction map  ${\mathbb C}[T]\to {\mathbb C}[Z(G^{s\circ})^\circ s]$ is surjective.
\end{enumerate}
\end{theorem}
\pf 1. The variety $Z(G^{s\circ})^\circ s/W(S)$ is the quotient of a smooth variety (a shifted torus) by the action of a finite group, hence it is normal. Every class in $\overline{J(s)}$ meets $T$ and $T\cap\overline{J(s)}=W\cdot (Z(G^{s\circ})^\circ s)$. Also, two elements in $T$ are $G$-conjugate if and only if they are $W$-conjugate, hence we have an isomorphism $ \overline{J(s)}//G\simeq W\cdot (Z(G^{s\circ})^\circ s)/W$ induced from the isomorphism $G//G\simeq T/W$. 

We consider the morphism $\gamma\colon Z(G^{s\circ})^\circ s/W(S)\to  W\cdot (Z(G^{s\circ})^\circ s)/W$ induced by $zs\mapsto W\cdot(zs)$. It is surjective by construction, bijective on the dense subset $(Z(G^{s\circ})^\circ s)^{reg}/W(S)$ and finite, since the intersection of $W\cdot (zs)$ with $Z(G^{s\circ})^\circ s$ is finite. Hence $\gamma$ is a normalisation morphism.

2. The variety $\overline{S}//G$ is normal if and only if the normalisation morphism is an isomorphism. This happens if and only if the composition  
\begin{equation*}
Z(G^{s\circ})^\circ s/W(S)\simeq \overline{S}//G\subseteq G//G\simeq T/W
\end{equation*}
is a closed embedding, i.e., if and only if  the corresponding algebra map between the rings of regular functions is surjective.
\hfill$\Box$

\section{An example: sheets and their quotients in  type $G_2$}

We list here the sheets in $G$ of type $G_2$ and all the conjugacy classes they contain. We shall denote by $\alpha$ and $\beta$, respectively, the short and the long simple roots. Since $G$ is adjoint, by \cite[Theorem 4.1]{gio} the sheets in $G$ are in bijection with $G$-conjugacy classes of pairs $(M,u)$ where $M$ is a pseudo-Levi subgroup of $G$ and $u$ is a rigid unipotent element in $M$. The corresponding sheet is $\overline{J(su)}^{reg}$ where $s$ is a semisimple element whose connected centralizer is $M$. The conjugacy classes of pseudo-Levi subgroups of $G$ are those corresponding to the following subsets $\Pi$ of the extended Dynkin diagram:
\begin{enumerate}
\item $\Pi=\emptyset$, so $M=T$, $u=1$, $s$ is a regular semisimple element and $S$ consists of all regular conjugacy classes;
\item $\Pi=\{\alpha\}$. Here $[M,M]$ is of type $\tilde{A}_1$, so $u=1$ and $s=\alpha^\vee(\zeta)\beta^\vee(t^2)=(3\alpha+2\beta)^\vee(\zeta^{-1})$ for $\zeta\neq0,\,\pm 1$;
\item $\Pi=\{\beta\}$. Here $[M,M]$ is of type ${A}_1$ so $u=1$ and $s=\alpha^\vee(\zeta^2)\beta^\vee(\zeta^3)=(2\alpha+\beta)^\vee(\zeta)$ for $\zeta\neq0,\,1\,e^{{2\pi i}/{3}},\,e^{-{2\pi i}/{3}}$;
\item $\Pi=\{\alpha_0,\beta\}$. Here $[M,M]$ is of type ${A}_2$ so $u=1$; the corresponding $s=(2\alpha+\beta)^\vee(e^{2\pi i/3})$ is isolated and $S=G\cdot s$;
\item $\Pi=\{\alpha_0,\alpha\}$. Here $[M,M]$ is of type $\tilde{A}_1\times A_1$ so $u=1$, the corresponding $s=(3\alpha+2\beta)^\vee(-1)$ is isolated and $S=G\cdot s$;
\item $\Pi=\{\alpha,\,\beta\}$ so $L=G$. In this case we have three possible choices for  $u$  rigid unipotent, namely 1, $x_{\alpha}(1)$ or $x_{\beta}(1)$ (cfr. \cite{spaltenstein}). Each of these classes is a sheet on its own.
\end{enumerate}

The only sheets containing more than one conjugacy classes  are the regular one $S_0=G^{reg}$ corresponding to $\Pi=\emptyset$ and the two subregular ones, corresponding to $\Pi_1=\{\alpha\}$ and $\Pi_2=\{\beta\}$. For $S_0$ we have $Z^\circ s=T$, $W(S)=W$ so $S_0/G$ is in bijection with $T/W$ and $\overline{S_0}//G\simeq G//G$ which is normal.  For $S_1$ and $S_2$ we have:
$$\begin{array}{ll}
&S_1=\overline{J((3\alpha+2\beta)^\vee(\zeta_0))}^{reg}\\
&=\left(\bigcup_{\zeta^2\neq0,1}G\cdot (3\alpha+2\beta)^\vee(\zeta)\right)\cup{\rm Ind}_{\tilde{A}_1}^G(1)\cup G\cdot \left((3\alpha+2\beta)^\vee(-1){\rm Ind}_{\tilde{A}_1}^{A_1\times\tilde{A}_1}(1)\right)\\
&=\left(\bigcup_{\zeta^2\neq0,1}G\cdot(3\alpha+2\beta)^\vee(\zeta)\right)\cup G\cdot\left((x_{\beta}(1)x_{\alpha_0}(1))\cup G\cdot (3\alpha+2\beta)^\vee(-1)x_{\alpha_0}(1)\right)
\end{array}$$ for $\zeta_0\neq0,\pm1$
and
$$\begin{array}{ll}
&S_2=\overline{J((2\alpha+\beta)^\vee(\xi_0))}^{reg}\\
&=\left(\bigcup_{\xi^3\neq0,1}G\cdot(2\alpha+\beta)^\vee(\xi)\right)\cup{\rm Ind}_{{A}_1}^G(1)\cup G\cdot \left((2\alpha+\beta)^\vee(e^{2\pi i/3}){\rm Ind}_{{A}_1}^{A_2}(1)\right)\\
&=\left(\bigcup_{\xi^3\neq0,1}G\cdot(2\alpha+\beta)^\vee(\xi)\right)\cup G\cdot(x_{\beta}(1)x_{\alpha_0}(1))\cup G\cdot \left((2\alpha+\beta)^\vee(e^{2\pi i/3}) x_{\alpha_0}(1)\right)
\end{array}$$
for some $\xi_0\neq0,1,e^{\pm 2\pi i/3}$.

In both cases $M$ is a Levi subgroup of a parabolic subgroup of $G$. By Lemmata \ref{lem:WS1} and \ref{WS} we have $W(S_1)=W(S_1)^u=\langle s_\alpha, s_{3\alpha+2\beta}\rangle$ and $W(S_2)=W(S_1)^u=\langle s_\beta, s_{2\alpha+\beta}\rangle$. Also $Z(M)^\circ=Z(M)^\circ s$  in both cases, so 
\begin{align*}
&S_1/G\simeq (3\alpha+2\beta)^\vee({\mathbb C^\times})/\langle s_\alpha, s_{3\alpha+2\beta}\rangle\simeq  (3\alpha+2\beta)^\vee({\mathbb C^\times})/\langle s_{3\alpha+2\beta}\rangle\\
&S_2/G\simeq (2\alpha+\beta)^\vee({\mathbb C^\times})/\langle s_\beta, s_{2\alpha+\beta}\rangle\simeq (2\alpha+\beta)^\vee({\mathbb C^\times})/\langle s_{2\alpha+\beta}\rangle,
\end{align*}
where the $\simeq$ symbols stand for the bijection $\overline{\theta}$.

Let us analyze normality of $\overline{S_1}//G$. Here, $Z(M)^\circ=(3\alpha+2\beta)^\vee({\mathbb C}^*)\simeq {\mathbb C}^*$, so ${\mathbb C}[Z(M)^\circ]^{W(S)}={\mathbb C}[\zeta+\zeta^{-1}]$.
On the other hand, since $G$ is simply connected, ${\mathbb C}[T]^W=( {\mathbb C}\Lambda)^W$ is the polynomial algebra generated by
$f_1=\sum_{\gamma\in\Phi\atop
\gamma\;\;{\rm  short }} e^\gamma$ and $f_2=\sum_{\gamma\in\Phi\atop
\gamma\;\;{\rm  long }} e^\gamma$, \cite[Ch.VI, \S 4, Th\'eor\`eme 1]{bourbaki}
Then, 
$$\rho(f_1)((3\alpha+2\beta)^\vee(\zeta))=f_1((3\alpha+2\beta)^\vee(\zeta))=\sum_{\gamma\in\Phi\atop
\gamma\;\;{\rm  short }} \zeta^{(\gamma, (3\alpha+2\beta)^\vee)}=2+2\zeta+2\zeta^{-1}$$ so the restriction map is surjective and $\overline{S_1}//G$ is normal.

Let us consider normality of $\overline{S_2}//G$. Here, $Z(M)^\circ=(2\alpha+\beta)^\vee({\mathbb C}^*)\simeq{\mathbb C}^*$, so ${\mathbb C}[Z]^\Gamma={\mathbb C}[\zeta+\zeta^{-1}]$. 
Then, 
$$\rho(f_1)(2\alpha+\beta)^\vee(\zeta)=f_1((2\alpha+\beta)^\vee(\zeta))=\sum_{\gamma\in\Phi\atop
\gamma\;\;{\rm  short }} \zeta^{(\gamma, (2\alpha+\beta)^\vee)}=\zeta^2+\zeta^{-2}+2(\zeta+\zeta^{-1})$$
whereas 
$$\rho(f_2)(2\alpha+\beta)^\vee(\zeta)=f_2((2\alpha+\beta)^\vee(\zeta))=\sum_{\gamma\in\Phi\atop
\gamma\;\;{\rm  long }} \zeta^{(\gamma, (2\alpha+\beta)^\vee)}=2+2\zeta^3+2\zeta^{-3}.$$
 Let us write $y=\zeta+\zeta^{-1}$. Then, $(\zeta^2+\zeta^{-2})=y^2-2$ and $\zeta^3+\zeta^{-3}=y^3-3y$ so ${\rm Im}(\rho)={\mathbb C}[y^2+2y,y^3-3y]={\mathbb C}[(y+1)^2,y^3+3y^2+6y+3-3y]={\mathbb C}[(y+1)^2,(y+1)^3]$.  
 Hence, $\rho$ is not surjective and $\overline{S_2}//G$ is not normal.
 
 We observe that  ${\rm Im}(\rho)$ is precisely the identification of the coordinate ring of $\overline{S_2}//G$ in ${\mathbb C}[T]^W$. We may thus see where this variety is not normal. We have:
${\rm Im}(\rho)={\mathbb C}[(y+1)^2,(y+1)^3]\cong {\mathbb C}[Y,Z]/(Y^3-Z^2)$ so this variety is not normal at $y+1=0$, that is, for $\zeta+\zeta^{-1}+1=0$. 
This corresponds precisely to the closed, isolated orbit $G\cdot  ((2\alpha+\beta)^\vee(e^{2\pi i/3}))x_{\alpha_0}(1)=G\cdot((2\alpha+\beta)^\vee(e^{-2\pi i/3}))x_{\alpha_0}(1)$. This example shows two phenomena: the first is that even if the sheet corresponsing to the set $\Pi_2$ in ${\rm Lie}(G)$ has a normal quotient \cite[Theorem 3.1]{broer}, the same does not hold in the group counterpart. The second phenomenon is that the non-normality locus corresponds to an isolated class in $\overline{S}_2$. In a forthcoming paper we will address the general problem of normality of $\overline{S}//G$ and we will prove and make use of the fact that  if the categorical quotient of the closure a sheet in $G$ is not normal, then it is certainly not normal at some isolated class.

\end{document}